\documentclass[12pt]{article}

\title{Nonparametric estimation for an autoregressive model}

\author{ Ouerdia Arkoun
\thanks{Laboratoire de Math\'ematiques Rapha\"el Salem,
UMR 6085 CNRS, Universit\'e de Rouen,
Avenue de l'Universit\'e, BP.12,
76801 Saint Etienne du Rouvray (France).\newline
email:  Ouerdia.Arkoun@univ-rouen.fr}
\and Serguei Pergamenchtchikov
\thanks{Laboratoire de Math\'ematiques Rapha\"el Salem,
UMR 6085 CNRS, Universit\'e de Rouen,
Avenue de l'Universit\'e, BP.12,
76801 Saint Etienne du Rouvray (France).\newline
email:  Serge.Pergamenchtchikov@univ-rouen.fr}
}

\usepackage{amssymb}
\usepackage{amsfonts}
\usepackage{amsmath}
\usepackage{amsthm}
\newtheorem{theorem}{Theorem}[section]

\newtheorem{lemma}[theorem]{Lemma}

\newtheorem{remark}[theorem]{Remark}

\newcommand\e{\varepsilon}
\def\bbr{{\mathbb R}}

\newcommand{\wh}{\widehat}
\newcommand{\wt}{\widetilde}

\newcommand\cH{{\cal H}}

\newcommand\cN{{\cal N}}
\newcommand\cP{{\cal P}}

\newcommand\cR{{\cal R}}

\newcommand\cU{{\cal U}}

\def\text#1{\hbox{#1}}
\def\proof{{\noindent \bf Proof. }}

\def\endproof{\mbox{\hfill $\qed$}}

\def\E{{\bf E}}
\def\P{{\bf P}}

\def\C{{\bf C}}

\def\Chi{{\bf 1}}
\def\d{\mbox{d}}
\def\build #1_#2{\mathrel{\mathop{\kern 0pt #1}\limits_{#2}}} %\textheight=9in
\newcommand{\zs}[1]{{\mathchoice{#1}{#1}{\lower.25ex\hbox{$\scriptstyle#1$}}
{\lower0.25ex\hbox{$\scriptscriptstyle#1$}}}}

\numberwithin{equation}{section}

\begin{document}

\maketitle

\begin{abstract}
 The paper deals with the nonparametric estimation problem 
at a given fixed point for an
autoregressive model with unknown distributed noise. Kernel estimate modifications are
proposed.
Asymptotic minimax and efficiency properties for  proposed estimators are shown.
\end{abstract}

\par
{\bf AMS (2000) Subject Classification :}
primary 62G07,62G08; secondary 62G20.

\section{Introduction}

We consider the following nonparametric autoregressive model
\begin{equation}\label{1.1}
y_k=S(x_k)y_{k-1}\, +\xi_k\,,
\quad 1\le k\le n\,, 
\end{equation}
where $S(\cdot)$ is an unknown $\bbr\to\bbr$ function, 
$x_k=k/n$, $y_0$ is a constant and  the noise random variables
$(\xi_k)_\zs{1\le k\le n}$ are i.i.d. with $\E \xi_k = 0$ and $\E \xi_k^2 = 1$.

The model \eqref{1.1} is a generalization of autoregressive processes of the first order.
In \cite{Da} the process \eqref{1.1} is considered with the function 
$S$ having a parametric form.
 Moreover, the paper \cite{Da-1} studies spectral properties
of the  stationary process \eqref{1.1} with the nonparametric
function $S$.

This paper deals with a nonparametric estimation of the autoregression 
coefficient function $S$ at a given point $z_0$, when the smoothness of $S$ is known. For this problem we make use of the following modified kernel estimator
\begin{equation}\label{1.2}
\hat{S}_\zs{n}(z_0)=\frac{1}{A_\zs{n}}\,
\sum^n_{k=1}\,Q(u_k)\,y_{k-1}\,y_\zs{k}\,
\Chi_\zs{(A_\zs{n}\ge d)},
\end{equation}
where $Q(\cdot)$ is a kernel function, 
$$
A_\zs{n}=\sum^n_{k=1}\,Q(u_k)y_{k-1}^2
\quad\mbox{with}\quad
u_k = \frac{x_k-z_0}{h}\,;
$$
 $d$ and $h$ are some positive parameters.

First we assume that the unknown function $S$ belongs to 
the {\em stable local H\"older class} at the point $z_\zs{0}$ with a known regularity 
$1\le \beta<2$.
This class will be defined below.
We find  an asymptotical (as $n\to\infty$) positive lower bound for the minimax risk
with the normalyzing coefficient 
\begin{equation}\label{1.3}
\varphi_\zs{n}=n^{\frac{\beta}{2\beta+1}}\,.
\end{equation}
To obtain this convergence rate we set
in \eqref{1.2}
\begin{equation}\label{1.4}
h=n^{-\frac{1}{2\beta+1}}
\quad\mbox{and}\quad
d=\kappa_\zs{n}\,nh\,,
\end{equation}
where  $\kappa_\zs{n}\ge 0$,
\begin{equation}\label{1.5}
\lim_\zs{n\to\infty}\kappa_\zs{n}=0
\quad\mbox{and}\quad
\lim_\zs{n\to\infty}\frac{h}{\kappa^2_\zs{n}}=0\,.
\end{equation}
As to the the kernel function we assume that
\begin{equation}\label{1.6}
\int_{-1}^{1}\, Q(z)\, \d z > 0 
\quad\mbox{and}\quad
 \int_{-1}^{1}\, z\,Q(z)\, \d z = 0 \,.
\end{equation}

In this paper we show that the estimator \eqref{1.2}
with the parameters \eqref{1.4}--\eqref{1.6}
 is asymptotically minimax, i.e. we show that the
asymptotical upper bound for the minimax risk 
with respect to the
stable local H\"older class is finite.

At the next step we study sharp asymptotic properties for
the minimax  estimators \eqref{1.2}. 
To this end
similarly to \cite{GaPe}
we introduce 
 the
{\em weak stable local H\"older class}. In this case we find
a positive constant giving
the exact asymptotic lower bound for the minimax risk with the
normalyzing coefficient \eqref{1.3}. Moreover, we show that
for the  estimator \eqref{1.2} with
the parameters \eqref{1.4}--\eqref{1.5} and
 the indicator kernel 
$Q=\Chi_\zs{[-1,1]}$ the asymptotic upper bound of the minimax
risk coincides with this constant, i.e. in this case such estimators
are asymptotically efficient.

In \cite{Be},  Belitser  consider the above  model with lipshitz condtions. The autor proposed a recursive estimator , and consider the estimation problem in a fixed t. By the quadratic risk, Belitser establish the convergence rate witout showing it's optimality.

Moulines and al  in \cite{Mou}, show that the convergence rate is optimal for the quadratic risk by using a recursive method for autoregressive model of order d.
We note that in our paper we establish an optimal convergence rate but the risk considered is different from the one used in \cite{Mou}, and assymptions are weaker then those of \cite{Mou}.

The paper is organized as follows. In the next section we give the main results.
In Section 3 we find asymptotical lowers bounds for the  minimax risks.
 Section 4 is devoted to  uppers bounds. Appendix contains some technical results.

\medskip

\section{ Main results}

Fisrt of all we assume that
the noise  in the model \eqref{1.1}, i.e. the i.i.d. random
variables  $(\xi_k)_\zs{1\le k\le n}$  have a density  $p$
(with respect to the Lebesue measure) 
 from the functional class $\cP$
defined as
\begin{align}\nonumber
\cP:=&
\left\{
p\ge 0\,:\,\int^{+\infty}_{-\infty}\,p(x)\,\d x=1\,,\quad
\int^{+\infty}_{-\infty}\,x\,p(x)\,\d x= 0 \,,\right.
\\[2mm] \label{2.1} 
&\quad\left.
\int^{+\infty}_{-\infty}\,x^2\,p(x)\,\d x= 1
\quad\mbox{and}\quad
 \int^{+\infty}_{-\infty}\,|x|^4\,p(x)\,\d x\le \sigma^*
\,
\right\}
\end{align}
with $\sigma^* \ge 3$.  Note that the $(0,1)$-gaussian density 
belongs to $\cP$. In the sequel we denote this density by $ p_\zs{0}$.

The problem is to  estimate the function $S(\cdot)$ 
at a fixed point $z_0\in ] 0, 1[$, i.e. the value $S(z_0)$.
For this problem we make use of the risk proposed in \cite{GaPe}. Namely,
  for any estimate $\wt{S}=\wt{S}_\zs{n}(z_0)$ 
(i.e. any mesurable with respect to the observations $(y_k)_\zs{1\le k\le n}$ function)
we set
\begin{equation}\label{2.2}
\cR_\zs{n}(\wt{S}_n,S)=
\sup_{p\in\cP} \,\E_\zs{S,p}|\wt{S}_\zs{n}(z_0)-S(z_0)|\,,
\end{equation}
where $\E_\zs{S,p}$ is the expectation taken with respect to the 
distribution $\P_\zs{S,p}$ of the vector $(y_1,...,y_n)$ in \eqref{1.1} corresponding to the function S and the density $p$ from $\cP$.

To obtain a stable (uniformly with respect to the function $S$ ) model \eqref{1.1}
 we assume (see \cite{Da} and \cite{Da-1}) 
that for some fixed $0<\e<1$
the unknown function $S$ belongs to the {\em stability set} 
\begin{equation}\label{2.3}
\Gamma_\zs{\e}=\{S\in \C_\zs{1}[0,1]\,:\, \|S\| \leq 1-\e\}\,,
\end{equation}
where $\|S\| = \sup_{0\leq x \leq1}|S(x)|$. Here
$\C_\zs{1}[0,1]$ is the Banah space of continuously differentiable
$[0,1]\to\bbr$ functions.

For fixed  constants $K>0$ and $0\le\alpha <1$  we define
the corresponding {\em stable local H\"older class} at the point $z_\zs{0}$
as
\begin{equation}\label{2.4}
 \cH^{(\beta)}(z_\zs{0},K,\e) 
= \left\{S \in \Gamma_\zs{\e} \,:\, \|\dot{S}\|\le K
\quad\mbox{and}\quad \Omega^*(z_\zs{0},S) \le K \right\}
\end{equation}
with $\beta=1+\alpha$ and
$$
\Omega^*(z_\zs{0},S) = \sup_{x\in [0,1]} 
\frac{|\dot{S}(x) - \dot{S}(z_\zs{0})|}{|x-z_\zs{0}|^\alpha}\,.
$$
First we show that the sequence \eqref{1.3} gives the optimal convergence rate
for the functions $S$ from $\cH^{(\beta)}(z_\zs{0},K,\e)$. 
We start with a lower bound.
\begin{theorem}\label{Th.2.1}
 For any  $K>0$ and $0<\e<1$
\begin{equation}\label{2.5}
\underline{\lim}_{n\to\infty} \, \inf_{\wt{S}}\quad
\sup_{S\in \cH^{(\beta)}(z_\zs{0},K,\e)}\, \varphi_n\, \cR_\zs{n}(\wt{S}_\zs{n},S) > 0,
\end{equation}
where the infimum is taken over all estimators.
\end{theorem}
\noindent Now we obtain an upper bound for the kernel estimator \eqref{1.2}.

\begin{theorem}\label{Th.2.2}
For any $K>0$ and $0<\e<1$
the kernel estimator \eqref{1.2} with the parameters
\eqref{1.4}--\eqref{1.6}  satisfies the following
inequality 
\begin{equation}\label{2.5-2}
\overline{\lim}_{n \to \infty} \, \sup_{S\in  \cH^{(\beta)}(z_\zs{0},K,\e)}\,
\varphi_n\,
\cR_\zs{n}(\hat{S}_\zs{n},S) < \infty.
\end{equation}
\end{theorem}
Theorem~\ref{Th.2.1} and  Theorem~\ref{Th.2.2} imply that the sequence
\eqref{1.3} is the optimal (minimax) convergence rate for any stable H\"older
class of regularity $\beta$, i.e. the estimator 
\eqref{1.2} with the parameters \eqref{1.4}--\eqref{1.6} is minimax
with respect to the functional class \eqref{2.4}.

Now we study some efficience properties 
for the minimax estimators \eqref{1.2}.
To this end similarly to \cite{GaPe} we make use  of the family of the {\em weak stable local H\"older classes} at the point $z_\zs{0}$, i.e. for any $\delta>0$  we set
\begin {equation}\label{2.6}
  \cU^{(\beta)}_{\delta,n}(z_\zs{0},\e) =
\left\{S\in \Gamma_\zs{\e}\,:\,
\|\dot S\|\le \delta^{-1}
\quad\mbox{and}\quad
|\Omega_\zs{h}(z_\zs{0},S)| \le \delta h^{\beta}
\right\}\,,
\end{equation}
where 
$$
\Omega_\zs{h}(z_\zs{0},S) =   \int_{-1}^{1}(S(z_0+uh)- S(z_0))\,\d u
$$
and $h$ is given in \eqref{1.4}.

\noindent Moreover, we set 
\begin{equation}\label{2.7}
\tau(S) = 1-S^2(z_0)\,.
\end{equation}
With the help of this function we 
describe the sharp lower bound for the minimax risks in this case.
\begin{theorem}\label{Th.2.3}
For any  $\delta>0$ and $0<\e<1$
 \begin{equation}\label{2.8}
\underline{\lim}_{n\to\infty} \, \inf_{\wt{S}}\,
\sup_{S\in  \cU^{(\beta)}_\zs{\delta,n}(z_\zs{0},\e)} \tau^{-1/2}(S)\,
 \varphi_n\,
\cR_\zs{n}(\wt{S}_\zs{n},S) \ge
 \E |\eta|\,,
\end{equation}
where $\eta$ is a gaussian random variable with the parameters $(0,1/2)$.
\end{theorem}

\begin{theorem}\label{Th.2.4}
The estimator \eqref{1.2} with the parameters \eqref{1.4}--\eqref{1.5} and $Q(z) = 1_{[-1,1]}$
 satisfies the following
inequality 
$$
\overline{\lim}_{\delta \to 0}\,
\overline{\lim}_{n\to\infty}\,
\sup_{S\in  \cU^{(\beta)}_\zs{\delta,n}(z_\zs{0},\e)} \tau^{-1/2}(S)\,
\varphi_n \,
\cR_\zs{n}({\hat{S}}_n,S) \le \E |\eta|\,,
 $$
where $\eta$ is a gaussian random variable with the parameters $(0,1/2)$. 
\end{theorem}
Theorems \ref{Th.2.3} and \ref{Th.2.4} imply that the  estimator \eqref{1.2}, \eqref{1.4}--\eqref{1.5}
 with the indicator
kernel is asymptotically efficient.
\begin{remark}\label{Re.2.1}
One can show (see \cite{GaPe}) that
for any $0<\delta<1$ and $n\ge 1$
$$
  \cH^{(\beta)}(z_\zs{0},\delta,\e)\subset  \cU^{(\beta)}_{\delta,n}(z_\zs{0},\e)\,.
$$ 
This means that the ``natural'' normalyzing coefficient for the functional class \eqref{2.6}
is the sequence \eqref{1.3}. Theorem~\ref{Th.2.3} and  Theorem~\ref{Th.2.4} extend usual the
H\"older approach for the point estimation  by keeping the minimax convergence rate \eqref{1.3}.
\end{remark}

\section{Lower bounds}
\subsection{Proof of Theorem~\ref{Th.2.1}}
 Note that to prove \eqref{2.5} it suffices to show that
\begin{equation}\label{3.1}
\underline{\lim}_{n\to\infty} \,\inf_{\wt{S}}
\sup_{S\in \cH^{(\beta)}(z_\zs{0},K,\e)}\,
\E_\zs{S,p_\zs{0}}\,\psi_n(\wt{S}_n,S)\,>0\,,
\end{equation}
where
$$
\psi_n(\wt{S}_n,S) = \varphi_n|\wt{S}_n(z_0)- S(z_0)|\,.
$$

We make use of the similar method proposed by Ibragimov and Hasminskii 
to obtain a lower bound for the density estimation problem in \cite{IbHa}.
First we chose the corresponding parametric family in $\cH^{(\beta)}(z_\zs{0},K,\e)$.
 Let $V$ be a two times continuously differentiable function
 such that 
$\int_{-1}^{1} V(z) dz > 0$ and $ V(z) = 0 $ for any $ |z| \ge 1$.
We set
\begin{equation}\label{3.2}
S_\zs{u}(x)= \frac{u}{\varphi_n}\,
V \left(\frac{x-z_0}{h}\right)\,,
\end{equation}
where  $\varphi_n$ and $h$ are defined in \eqref{1.3} and \eqref{1.4}.

It is easy to see that for any
$z_0-h\le x\le z_0+h$ 
\begin{align*}
 |\dot{S}_u(x) - \dot{S}_u(z_\zs{0})|
& =
\frac{|u|}{h\varphi_n} \left|\dot{V}\left( \frac{x-z_0}{h}\right) - 
\dot{V}(0)\right| \\[2mm]
& \leq \frac{|u|}{h\varphi_n}\, V''_\zs{*}\,\left|\frac{x-z_0}{h}\right|\,
 \le\, |u| V''_\zs{*}\,|x-z_0|^\alpha\,,
\end{align*}
where $V''_\zs{*}=\max_{|z|\leq 1}|\ddot{V}(z)|$. Therefore, 
for all $0 < u\le u^* = K/V''_\zs{*}$
we obtain that
$$
\sup_\zs{z_0-h\le x\le z_0+h}
\frac{|\dot{S}_u(x) - \dot{S}_u(z_\zs{0})|}{|x-z_0|^\alpha}\,\le \,K\,.
$$
Moreover, by the definition \eqref{3.2} for all $x> z_\zs{0}+h$ 
$$
\dot{S}_\zs{u}(x)=\dot{S}_\zs{u}(z_\zs{0}+h)=0
\quad\mbox{and}\quad
\dot{S}_\zs{u}(x)=\dot{S}_\zs{u}(z_\zs{0}-h)=0
$$
for all $x< z_\zs{0}-h$ respectively. Therefore, the last inequality implies  that
$$
\sup_\zs{|u|\le u^*}
\Omega^*(z_\zs{0},S_\zs{u}) \le K\,,
$$
where the function $\Omega^*(z_\zs{0},S)$ is defined in \eqref{2.4}.

This means that there exists $n_\zs{K,\e}>0$ such that 
$S_\zs{u}\in \cH^{(\beta)}(z_\zs{0},K,\e)$ for all $|u|\le u^*$ and $n\ge n_\zs{K,\e}$. 
Therefore, for all 
$n\ge n_\zs{K,\e}$
and for any estimator $\wt{S}_\zs{n}$ we estimate with below the supremum
 in \eqref{3.1} as
\begin{align}\nonumber
\sup_{S\in \cH^{(\beta)}(z_\zs{0},K,\e)}\,\E_\zs{S,p_\zs{0}}\,
\psi_n(\wt{S}_n,S)&
\ge \sup_{|u|\le u^*}\E_\zs{S_\zs{u},p_\zs{0}}\psi_\zs{n}(\wt{S}_\zs{n},S_\zs{u})\\
\label{3.3}
&\ge
\frac{1}{2 b}\int^{b}_{-b}\E_\zs{S_\zs{u},p_\zs{0}}\psi_\zs{n}(\wt{S}_\zs{n},S_\zs{u})\d u 
\end{align}
for any $0<b\le u^*$.

Notice that for any $S$
the measure $\P_{S,p_0}$ is equivalent to the measure $\P_\zs{0,p_\zs{0}}$, where $\P_{0,p_0}$ is the distribution of the vector $(y_1,\ldots,y_n)$ in \eqref{1.1} 
corresponding to the function $S=0$
and the gaussian $(0,1)$ noise density $p_\zs{0}$, i.e. the random variables 
$(y_1,...,y_n)$ are i.i.d. $\cN(0,1)$ with respect to the measure $\P_\zs{0,p_\zs{0}}$.
In the sequal we denote $\P_{0,p_0}$ by $\P$.  
It is easy to see that in this case the Radom-Nikodym derivative can be written as
$$
\rho_n(u)=\frac{\d \P_{S_\zs{u},p_0}}{\d \P}
=
e^{u\,\varsigma_n\,\eta_n
-\frac{u^2}{2}\varsigma^2_n}\
$$
with
$$
\varsigma^2_n=\frac{1}{\varphi^2_n}
\sum^n_{k=1}\,V^2(u_k) \xi^2_{k-1} 
\quad\mbox{and}\quad
\eta_n=\frac{1}{\varsigma_n\,\varphi_n}\,
\sum^n_{k=1}\,V(u_k)\,\xi_{k-1}\,\xi_k
\,.
$$
Through the large numbers law we obtain 
$$
\P-\lim_{n\to\infty}\,\varsigma^2_n
 = \lim_{n\to\infty}
\frac{1}{nh}
\sum^{k^*}_{k=k_*}\, V^2\,(u_k) \xi^2_{k-1}
 =\int_{-1}^1 V^2(u) du = \sigma^2\,,
$$
where 
\begin{equation}\label{3.3-1}
k_*= [nz_0 - nh ] + 1
\quad\mbox{and}\quad
k^*=  [nz_0 + nh ]\,.
\end{equation}
Here $[a]$ is the integer part of $a$.

Moreover, by the central limit theorem for martingales (see \cite{He} and \cite{Re}),
 it is easy to see that under the measure $\P$ 
$$
 \eta_n\quad \Longrightarrow\quad \cN(0,1)
\quad\mbox{as}\quad n \to \infty
\,.
$$
Therefore we represent the Radon-Nykodim density in the following asymptotic form
$$
\rho_n(u)
=e^{u\sigma\eta_n-\frac{u^2\sigma^2}{2}+r_n}\,,
$$
where 
$$
\P - \lim_{n\to\infty} r_n=0\,.
$$
This means that in this case the Radon-Nikodym density 
$(\rho_n(u))_\zs{n\ge 1}$ satisfies the L.A.N. property
and we can make use the method from theorem 12.1 of \cite{IbHa} to obtain the folowing inequality
 \begin{equation}\label{3.4}
\underline{\lim}_{n\to\infty}\inf_{\wt{S}} \,
\frac{1}{2 b}\int^{b}_{-b}\E_\zs{S_\zs{u},p_\zs{0}}\psi_\zs{n}(\wt{S}_\zs{n},S_\zs{u})\,\d u 
\,\ge\,  I(b,\sigma)\,,
 \end{equation}
where
$$
I(b,\sigma)  = \frac{\max (1,b-\sqrt{b})}{b}\frac{\sigma}{\sqrt{2\pi}}\int_{-\sqrt{b}}^{\sqrt{b}} e^{-\sigma^2\frac{u^2}{2}} du
$$
and $0<b\le u^*$.
Therefore, inequalities \eqref{3.3} and \eqref{3.4} imply \eqref{3.1}. 
Hence Theorem~\ref{Th.2.1}.
\endproof

\subsection{Proof of Theorem \ref{Th.2.3}}

First, similarly to the proof of Theorem~\ref{Th.2.1} we choose the 
corresponding parametric functional
family $S_\zs{u,\nu}(\cdot)$ in the form \eqref{3.2} with the function $V=V_\zs{\nu}$ defined as
$$
V_\nu(x)=\nu^{-1}\int^\infty_{-\infty}\wt{Q}_\nu(u) g\left(
\frac{u-x}{\nu}\right) \d u\,,
$$
 where $\wt{Q}_\nu(u)={\bf 1}_{\{|u|\le 1-2\nu\}}+
2{\bf 1}_{\{1-2\nu\le |u|\le 1-\nu\}}$ with $0<\nu<1/4$ and $g$ is some even nonnegative infinitely differentiable
function such that $g(z)=0$ for $|z|\ge 1$ and $\int^1_{-1}\,g(z)\ \d z=1$. 
 One can show (see \cite{GaPe}) that for any $b>0$, $0<\delta<1$
and $0<\nu<1/4$ 
there exists $n_\zs{*}=n_\zs{*}(b,\delta,\nu)>0$ such that for all
$|u|\le b$ and $n\ge n_\zs{*}$
$$
S_\zs{u,\nu} \in 
\cU^{(\beta)}_{\delta,n}(z_\zs{0},\e)\,.
$$
Therefore, in this case for any $n\ge n_\zs{*}$
\begin{align*}
\varphi_n\,
\sup_{S\in  \cU^{(\beta)}_\zs{\delta,n}(z_\zs{0},\e)} \tau^{-1/2}(S)\,\cR_\zs{n}(\wt{S}_\zs{n},S)
&\ge 
\sup_{S\in  \cU^{(\beta)}_\zs{\delta,n}(z_\zs{0},\e)} \tau^{-1/2}(S)
\E_\zs{S,p_\zs{0}}\,\psi_n(\wt{S}_n,S)\\[2mm]
&
\ge\,
\tau_\zs{*}(n,b)\,
\frac{1}{2 b}\int^{b}_{-b}\,
\E_\zs{S_\zs{u,\nu},p_\zs{0}}\psi_\zs{n}(\wt{S}_\zs{n},S_\zs{u,\nu})\d u \,.
\end{align*}
where
$$
\tau_\zs{*}(n,b)=
\inf_\zs{|u|\le b}\,\tau^{-1/2}(S_\zs{u,\nu})\,.
$$
The definitions \eqref{2.7} and \eqref{3.2} imply
that for any $b>0$
$$
\lim_\zs{n\to\infty}\sup_\zs{|u|\le b}\,
|\tau(S_\zs{u,\nu})-1|\,=\,0\,.
$$
Therefore, by the same way as in the proof of Theorem~\ref{Th.2.1} we obtain that for any 
$b>0$ and $0<\nu<1/4$
\begin{equation}\label{3.5}
\underline{\lim}_{n\to\infty} \,\inf_{\wt{S}}\,
\sup_{S\in  \cU^{(\beta)}_\zs{\delta,n}(z_\zs{0},\e)} \tau^{-1/2}(S)\,
 \varphi_n\,
\cR_\zs{n}(\wt{S}_\zs{n},S)
\ge I(b,\sigma_\zs{\nu})\,,
\end{equation}
where the function $I(b,\sigma_\zs{\nu})$ is defined in \eqref{3.4} 
with $\sigma^2_{\zs\nu}=\int_{-1}^{1}\,V^2_\nu(u)\,\d u$. It is easy to check that
 $\sigma^2_{\zs\nu} \rightarrow 2 $ as $\nu \rightarrow 0$.
Limiting $b\to\infty$  and $\nu \rightarrow 0$ in \eqref{3.5} yield
the inequality \eqref{2.8}.
Hence Theorem~\ref{Th.2.3}.

\endproof

\section{Upper bounds}
\subsection{Proof of Theorem \ref{Th.2.2}}
First of all we set
\begin{equation}\label{4.1}
\wt{A}_\zs{n}=\frac{A_\zs{n}}{\varphi^2_\zs{n}}
\quad\mbox{and}\quad
\wh{A}_\zs{n}=
\frac{1}{\wt{A}_\zs{n}}\,\Chi_\zs{(\wt{A}_\zs{n}> \kappa_\zs{n})}\,.
\end{equation}
Now from \eqref{1.2} we represent the estimate error as
\begin{equation}\label{4.2}
\hat{S}_\zs{n}(z_0)-S(z_0)=
-S(z_0)\,\Chi_\zs{(\wt{A}_\zs{n}\le \kappa_\zs{n})}
+
\frac{1}{\varphi_n}\,
\wh{A}_\zs{n}\,\zeta_\zs{n}
+
\frac{1}{\varphi_n}\,
\wh{A}_\zs{n}\,B_\zs{n}\,,
\end{equation}
with
$$
\zeta_\zs{n}=\frac{\sum^n_{k=1}\,Q(u_k)y_{k-1}\,\xi_k}{\varphi_n}
\quad\mbox{and}\quad
B_\zs{n}=\frac{\sum^n_{k=1}\,Q(u_k)\,(S(x_k)-S(z_0))y_{k-1}^2}{\varphi_\zs{n}}
\,.
$$
Note that, the first  term in 
the right hand of \eqref{4.2} is studied in Lemma~\ref{Le.A.3}.
To estimate the second term we make use of Lemma~\ref{Le.A.2} which 
implies directly
$$
\overline{\lim}_{n\to\infty}\,
\sup_{S\in\cH^{(\beta)}(z_\zs{0},K,\e)}\,
\sup_\zs{p\in\cP}\,
\E_\zs{S,p}\,\zeta^2_\zs{n}\,< \infty 
$$
and, therefore, by \eqref{A.8} we obtain 
$$
\overline{\lim}_{n\to\infty}\,
\sup_{S\in\cH^{(\beta)}(z_\zs{0},K,\e)}\,
\sup_\zs{p\in\cP}\,
\E_\zs{S,p}\,|\wh{A}_\zs{n}|\,|\zeta_\zs{n}|\,<\,\infty\,.
$$
Let us estimate now the last term in the right hand of \eqref{4.2}. To this end we need to show
that
\begin{equation}\label{4.3}
\overline{\lim}_{n\to\infty}\,
\sup_{S\in\cH^{(\beta)}(z_\zs{0},K,\e)}\,
\sup_\zs{p\in\cP}\,
\E_\zs{S,p}\,B^2_\zs{n}\,<\,\infty\,.
\end{equation}
Indeed, putting 
$ r_k = S(x_k)- S(z_0) - \dot{S}(z_0)(x_k-z_0)$
 by the Taylor Formula we represent $B_\zs{n}$ as
$$
B_\zs{n} 
= \frac{h}{\varphi_n}\dot{S}(z_0) \wt{B}_\zs{n} +\frac{1}{\varphi_n} \wh{B}_n\,,
$$
where 
$\wt{B}_\zs{n}=\sum^n_{k=1}\,Q(u_k)\,u_\zs{k}\, y^2_{k-1}$
and
 $\wh{B}_n=\sum^n_{k=1}\,Q(u_k)\,r_k \, y^2_{k-1}$.
We remind that by the condition \eqref{1.6}
$\int^1_\zs{-1}uQ(u)\d u=0$. Therefore through Lemma~\ref{Le.A.2} we obtain
$$
\lim_{n\to\infty}\,
\frac{h^2}{\varphi^2_n}\,
\sup_{S\in\cH^{(\beta)}(z_\zs{0},K,\e)}\,
\sup_\zs{p\in\cP}\,
\E_\zs{S,p}\,\wt{B}^2_\zs{n}\,=0\,.
$$
Moreover, for any function $S\in\cH^{(\beta)}(z_\zs{0},K,\e)$
and for $k_\zs{*}\le k\le k^*$ ($k_\zs{*}$ and $ k^*$ are given in \eqref{3.3-1})
$$
 |r_k| \,=\,\left| \int_{z_0}^{x_k}\left(\dot{S}(u) - \dot{S}(z_0)\right)\d u\right| \, \le K|x_k-z_0|^\beta 
\le K h^\beta=K\varphi^{-1}_\zs{n}\,,
$$
i.e. $\wh{B}_n\le \varphi_n\wt{A}_\zs{n}$. Therefore, by Lemma~\ref{Le.A.2}
$$
\overline{\lim}_{n\to\infty}\,
\sup_{S\in\cH^{(\beta)}(z_\zs{0},K,\e)}\,
\sup_\zs{p\in\cP}\, \frac{1}{\varphi^2_n}
\E_\zs{S,p}\,\wh{B}^2_\zs{n}\,<\infty\,.
$$
This implies \eqref{4.3}. Hence Theorem~\ref{Th.2.2}.
\endproof  

\subsection{Proof of Theorem \ref{Th.2.4}}
Similarly to  Lemma A.2 from \cite{GaPe} 
by making use of Lemma~\ref{Le.A.1} and Lemma~\ref{Le.A.2}
we can show that
$$
\sqrt{\frac{\tau(S)}{2}}\,\zeta_n\quad \Longrightarrow \quad \cN(0,1)
\quad\mbox{as}\quad
n\to\infty
$$
uniformly 
 in $S \in \Gamma_\zs{\e}$ and  $p \in \cP$. Therefore, by Lemma~\ref{Le.A.2} we obtain that
uniformly 
 in $S \in \Gamma_\zs{\e}$ and  $p \in \cP$
$$
\tau^{-1/2}(S)\,
\wh{A}_\zs{n}\,\zeta_\zs{n}\quad \Longrightarrow \quad
\cN\left(0\,,\,1/2\right)
\quad\mbox{as}\quad
n\to\infty
\,.
$$
Moreover, by applying the Burkh\"older inequality and Lemma~\ref{Le.A.2} to the martinagale $\zeta_\zs{n}$ 
we deduce that
$$
\overline{\lim}_{n\to\infty}\,
\sup_{S\in\cH^{(\beta)}(z_\zs{0},K,\e)}\,
\sup_\zs{p\in\cP}\,
\E_\zs{S,p}\,\zeta^4_\zs{n}\,< \infty \,.
$$
Therefore, inequality \eqref{A.8} implies that the sequence $(\wh{A}_\zs{n}\,\zeta_\zs{n})_\zs{n\ge 1}$
is uniformly integrable. This means that
$$
\lim_{n\to\infty}\,
\sup_{S\in\cH^{(\beta)}(z_\zs{0},K,\e)}\,
\sup_\zs{p\in\cP}\,
\left|\tau^{-1/2}(S)\,
\E_\zs{S,p}\,|\wh{A}_\zs{n}\,\zeta_\zs{n}|
-
\E|\eta|
\right|=0\,,
$$
where
$\eta$ is a gaussian random variable with the parameters $(0,1/2)$. 
Now to finish this proof we have to show that 
\begin{equation}\label{4.4}
\lim_\zs{\delta\to 0}\,
\overline{\lim}_{n\to\infty}\,
\sup_{S\in  \cU^{(\beta)}_\zs{\delta,n}(z_\zs{0},\e)}\,
\sup_\zs{p\in\cP}\,
\E_\zs{S,p}\,B^2_\zs{n}\,=\,0\,.
\end{equation}
Indeed, by setting 
$f_\zs{S}(u)= S(z_0+h u)-S(z_0)$
we rewrite $B_\zs{n}$ as
\begin{equation}\label{4.5}
B_\zs{n}
 = \frac{1}{\varphi_\zs{n}}\,\sum^{k^*}_{k=k_*}\,f_\zs{S}(u_k)\,y_{k-1}^2=
\varphi_\zs{n}\,\varrho_\zs{n}(f_\zs{S},S)\,+
\,\frac{\varphi_\zs{n}}{\tau(S)}\,
\Omega_\zs{h}(z_\zs{0},S)
\,,
\end{equation}
where
$$
\varrho_\zs{n}(f,S)=
\frac{\sum^{n}_{k=1}\,f(u_k)y^2_{k-1}}{\varphi^2_\zs{n}}\,-\,\frac{1}{\tau(S)}
\,\int_{-1}^{1}\, f(u)\d u
$$
and 
$\Omega_\zs{h}(z_\zs{0},S)$ is defined in
\eqref{2.6}. The definition \eqref{2.7} implies that for any $S\in\Gamma_\zs{\e}$
\begin{equation}\label{4.6}
\e^2\le \tau(S)\le 1
\,.
\end{equation}
From here by the definition \eqref{2.6} we obtain that
$$
|B_\zs{n}|\le \varphi_\zs{n}\,|\varrho_\zs{n}(f_\zs{S},S)|+
\frac{\delta}{\e^2}\,.
$$
Moreover, for any
$S\in\cU^{(\beta)}_\zs{\delta,n}(z_\zs{0},\e)$
the function $f_\zs{S}$ satisfies the following inequality
$$
\|f_\zs{S}\|+\|\dot{f}_\zs{S}\|\le \delta^{-1}\,h\,.
$$
We note also that $\varphi_\zs{n}h^2\to 0$ as $n\to\infty$.
Therefore, by making use of Lemma~\ref{Le.A.2} with $R=h/\delta$ we obtain \eqref{4.4}.
Hence Theorem~\ref{Th.2.4}.
\endproof

\renewcommand{\theequation}{A.\arabic{equation}}
\renewcommand{\thetheorem}{A.\arabic{theorem}}
\renewcommand{\thesubsection}{A.\arabic{subsection}}
\section{Appendix}

In this section we study  distribution properties of the stationary process \eqref{1.1}. 
\begin{lemma}\label{Le.A.1}
For any $0<\e<1$ the random variables \eqref{1.1} satisfy the following
moment inequality
\begin{equation}\label{A.1}
m^*=
\sup_\zs{n\ge 1}\,\sup_\zs{0\le k\le n}\,
\sup_\zs{S\in\Gamma_\zs{\e}}\,\sup_\zs{p\in\cP}\,
\E_\zs{S,p}\,y^4_\zs{k}\,<\,\infty\,.
\end{equation}
\end{lemma}
\proof 
One can deduce from \eqref{1.1} with $S\in\Gamma_\zs{\e}$
 that for all $1\le k\le n$
$$
y^4_\zs{k}\le \left((1-\e)^k|y_\zs{0}|+\sum^k_\zs{j=1}\,(1-\e)^{k-j}\,|\xi_\zs{j}|\right)^4
\le 8y^4_\zs{0}+
 8\left(
\sum^k_\zs{j=1}\,(1-\e)^{k-j}\,|\xi_\zs{j}|\right)^4
\,.
$$
Moreover, by the H\"older inequality with $q=4/3$ and $p=4$
$$
y^4_\zs{k}\,\le 8|y_\zs{0}|^4+\frac{8}{\e^3}\,
\sum^k_\zs{j=1}(1-\e)^{k-j}\,\xi^{4}_\zs{j}\,.
$$
Therefore, for any $p\in\cP$
$$
\E_\zs{S,p}\,y^4_\zs{k}\,\le\, 8\,|y_\zs{0}|^4\,+\,
\frac{8}{\e^4}\,\sigma_\zs{*}\,.
$$
Hence Lemma~\ref{Le.A.1}. \endproof

\noindent Now for any $K>0$ and $0<\e<1$ we set 
\begin{equation}\label{A.2}
\Theta_{K,\e} = \{ S\in \Gamma_\zs{\e}\,:\, \|\dot{S}\|\le K\}\,.
\end{equation}
\begin{lemma}\label{Le.A.2}
Let the function $f$ is two times continuously differentiable in $[-1,1]$, 
such that $f(u) =0$ for  $|u|\ge 1$. Then 
\begin{equation}\label{A.3} 
\overline{\lim_{n\to \infty}}\,\sup_\zs{R>0}
\frac{1}{(R h)^2}\sup_{\|f\|_\zs{1}\le R}\,
 \sup_{S \in \Theta_{K,\e}}\,\sup_{p\in\cP} \E_\zs{S,p}\, 
\varrho^2_\zs{n}(f,S)
 \,<\,\infty\,,
\end{equation}
where
$\|f\|_\zs{1} = \|f\| + \|\dot{f}\|$ and
$\varrho_\zs{n}(f,S)$
is defined in \eqref{4.5}.
\end{lemma}
\proof
First of all, note that
\begin{equation}\label{A.4} 
\sum^{n}_{k=1}\,f(u_k)y^2_{k-1}=
T_\zs{n}+a_\zs{n}\,,
\end{equation}
where
$$
T_\zs{n}=\sum^{k^*}_{k=k_*}\,f(u_k) y^2_\zs{k}
\quad\mbox{and}\quad
a_\zs{n}=\sum^{k^*}_\zs{k=k_\zs{*}}\,(f(u_\zs{k})-f(u_\zs{k-1}))\,y^2_\zs{k-1}
-f(u_\zs{k^*})\,y^2_\zs{k^*}
$$
with $k^*$ and $k_\zs{*}$ defined in \eqref{3.3-1}. Moreover, from the model \eqref{1.1} we find 
$$
T_\zs{n}= I_\zs{n}(f)+ \sum^{k^*}_{k=k_*}\,f(u_k)S^2(x_k)y^2_{k-1} + M_n\,,
$$
where
 $$
I_\zs{n}(f)=\sum^{k^*}_{k=k_*}\,f(u_k)
\quad\mbox{and}\quad
M_n  =\sum^{k^*}_{k=k_*}\,f(u_k)\,
\left(
2\,S(x_k)\,y_{k-1}\,\xi_k\,+\,\eta_k
\right)
$$
with $\eta_k = \xi_k^2- 1$.
By setting 
$$
C_n = \sum^{k^*}_{k=k_*}\, (S^2(x_k)-S^2(z_0))\,f(u_k)\,y^2_{k-1}
\quad\mbox{and}\quad
D_n = \sum^{k^*}_{k=k_*}\,f(u_k)(y^2_{k-1}-y^2_k)
$$ 
we get 
\begin{equation}\label{A.5} 
\frac{1}{\varphi^2_\zs{n}}\,
T_\zs{n}
= \frac{1}{\tau(S)}\, \frac{I_\zs{n}(f)}{\varphi^2_\zs{n}}
+ \frac{1}{\tau(S)} \frac{\Delta_\zs{n}}{\varphi^2_\zs{n}} 
\end{equation}
with $ \Delta_\zs{n} = M_\zs{n}+C_\zs{n} + S^2(z_0)\,D_\zs{n}$.
Moreover, taking into account that $\varphi^2_\zs{n}=nh$ we obtain 
\begin{align*} 
\frac{I_\zs{n}(f)}{\varphi^2_\zs{n}} &= \int_{-1}^{1}f(t)\d t
 +\sum^{k^*}_{k=k_*}\,\int_{u_{k-1}}^{u_k}f(u_k)\,\d t\,-\, \int_{-1}^{1}f(t)dt \\
& = 
\sum^{k^*}_{k=k_*}\,\int_{u_{k-1}}^{u_k} (f(u_k)- f(t)) dt + \int_{u_{k_{*}-1}}^{u_{k^*}} f(t)dt - \int_{-1}^{1}f(t)dt\,.
\end{align*}
We remind that $\|f\|+\|\dot{f}\|\le R$. Therefore 
$$
\left|\frac{1}{nh} \sum^{k^*}_{k=k_*}\,f(u_k)- \int_{-1}^{1}f(t)dt\right| 
 \le \frac{R}{nh}\,.
$$
Taking this into account in \eqref{A.5} and the lower bound for $\tau(S)$ given in \eqref{4.6} we 
find that
\begin{equation}\label{A.6}
\left|\frac{T_\zs{n}}{\varphi^2_\zs{n}}
-\frac{1}{\tau(S)} \int_{-1}^{1}f(t)dt\right| 
\le 
\frac{1}{\e^2}
\left(
\frac{R}{nh}+\frac{M_\zs{n}}{nh}+
\frac{C_\zs{n}}{nh}+\frac{D_\zs{n}}{nh}
\right)\,.
\end{equation}
Note that the sequence $(M_\zs{n})_\zs{n\ge 1}$ is a square integrable martingale.
Therefore,
\begin{align*}
\E_\zs{S,p}\,\left(\frac{1}{n h}M_n\right)^2
& = \frac{1}{(n h)^2}\,\E_\zs{S,p}\,
\sum^{k^*}_{k=k_*}\,f^2(u_k)\,
\left(
2\,S(x_k)\,y_{k-1}\,\xi_k\,+\,\eta_k
\right)^2\\
& 
\le \frac{4 R^2 \left(4\sqrt{m^*}+\sigma^*\right)}{n h},
\end{align*}
 where $m^*$ is given in \eqref{A.1}. 
Moreover, taking into account that $|S(x_k)-S(z_0)| \le L|x_k-z_0|$  
for any $S\in \Theta_\zs{L,\e}$ and that $k^*-k_\zs{*}\le 2n h$
we obtain that
\begin{align*}
 \frac{1}{(n h)^2}\,\E_\zs{S,p}\, C^2_n
& \le \frac{2}{n h}\sum^{k^*}_{k=k_*}\, |(S^2(x_k)-S^2(z_0))|^2\,f^2(u_k)\,\E_\zs{S,p}\,y^4_{k-1} \\
& \le \,16\,R^2\, L^2\, m^*\,h^2\,.
\end{align*}
Let us consider now the last term in the right hand of the inequality \eqref{A.6}. To this end
we make use  of the integration by parts formula, i.e. we represent $D_\zs{n}$ as
$$ 
D_n = \sum^{k^*}_{k=k_*}\,\left((f(u_k)-f(u_{k-1})\right)\,y^2_{k-1} +f(u_{k_{*}-1})\,y^2_{k_{*}-1} - f(u_{k^*})\,y^2_{k^*}\,.
$$
Therefore, taking into account that $\|f\|+\|\dot{f}\|\le R$
we obtain that

\begin{align*}
\E_\zs{S,p}\,D^2_\zs{n}
\le\, 3R^2\,  
\E_\zs{S,p}\,
\left(\frac{2}{n h}
\,
\sum^{k^*}_{k=k_*}\,y^4_{k-1}\,
+y^4_{k^*}+y^4_{k_{*}-1}
\right)
 \le 18\, R^2\,m^*\,.
\end{align*} 
By the same way we estimate the second term in the right hand of \eqref{A.4}.
Hence Lemma~\ref{Le.A.2}. 

\endproof

\begin{lemma}\label{Le.A.3}
The sequences $(\wt{A}_\zs{n})_\zs{n\ge 1}$ 
and $(\wh{A}_\zs{n})_\zs{n\ge 1}$ defined in \eqref{4.1}
 satisfy the following
properties
\begin{equation}\label{A.7}
\overline{\lim_{n\to \infty}}\,\frac{1}{h^2}
 \sup_{S \in \Theta_{K,\e}}\,\sup_{p\in\cP}\,
 \P_\zs{S,p}\,(\wt{A}_\zs{n}\le \kappa_\zs{n})
\,<\,\infty
\end{equation}
and
\begin{equation}\label{A.8}
\overline{\lim_{n\to \infty}}\,
 \sup_{S \in \Theta_{K,\e}}\,\sup_{p\in\cP}\,
 \E_\zs{S,p}\, 
\wh{A}^{4}_\zs{n}\,<\,\infty\,.
\end{equation}
\end{lemma}
\proof
It is easy to see that the inequality \eqref{A.7} follows directly from Lemma~\ref{Le.A.2}.
We check now the inequality \eqref{A.8}. By setting 
$\gamma_\zs{*}=\e^{-2}\int^1_\zs{-1}Q(u)\d u$ 
we get 
\begin{align*}
 \E_\zs{S,p}\, 
\wh{A}^{4}_\zs{n}\,&=4
\int^{\infty}_\zs{0}\,t^3\,
\P_\zs{S,p}\,\left(\wt{A}_\zs{n}\le t^{-1}\,,\,\wt{A}_\zs{n}> \kappa_\zs{n}\right)\,\d t\\
&\le 
4
\int^{\kappa^{-1}_\zs{n}}_\zs{0}\,t^3\,
\P_\zs{S,p}\,\left(\varrho_\zs{n}(Q,S)+\gamma_\zs{*}\le t^{-1}\right)\,\d t\\
&\le \left(\frac{2}{\gamma_\zs{*}}\right)^4+
\frac{1}{\kappa^4_\zs{n}}\,
\P_\zs{S,p}\,\left(|\varrho_\zs{n}(Q,S)|\ge \gamma_\zs{*}/2\right)\,.
\end{align*}
By making use of Lemma~\ref{Le.A.2} with
the condition \eqref{1.5} we obtain the inequality \eqref{A.8}.
\endproof

%%%%%%%%%%%%%%%%%%%%%%%%%%%%%%%%%%%%%%%%%%%%%%%%%%%%%%%%%%%%%%

\end{document}